\documentclass[11pt]{article}
\usepackage{amsmath}

\newtheorem{theorem}{Theorem}
\newtheorem{lemma}{Lemma}

\begin{document}
\begin{center}
{\bf TRANSVERSALS IN TREES}\\
\bigskip
Victor Campos\footnotemark[1]$^,$\footnotemark[3]\\
Va\v{s}ek Chv\'{a}tal\footnotemark[2]$^,$\footnotemark[4]\\
Luc Devroye\footnotemark[1]$^,$\footnotemark[5]\\
Perouz Taslakian\footnotemark[1]
\footnotetext[1]{School of Computer Science, McGill University,
  Montr\' eal, Qu\' ebec, Canada}
\footnotetext[2]{Canada Research Chair in Combinatorial Optimization,
Department of Computer Science and Software Engineering,
Concordia University, Montr\' eal, Qu\' ebec, Canada}
\footnotetext[3]{Supported by Conselho Nacional de Desenvolvimento
Cient\'\i fico e Tecnol\'{o}gico (CNPq), Brazil}
\footnotetext[4]{Supported by the Canada Research Chairs Program and by the
Natural Sciences and Engineering Research Council of Canada}
\footnotetext[5]{Supported by the
Natural Sciences and Engineering Research Council of Canada}
\end{center}

\bigskip

A {\em transversal\/} in a rooted tree is any set of nodes that meets
every path from the root to a leaf. We let $c(T,k)$ denote the number
of transversals of size $k$ in a rooted tree $T$. If $T$ has $n$ nodes
and $n\ge 2$, then
\begin{eqnarray*}
c(T,n) &=& 1,\\
c(T,n-1) &=& n,
\end{eqnarray*}
and
\begin{equation}\label{sandwich}
{\textstyle \binom{n-1}{k-1}} \;\le\; c(T,k) \;\le\; {\textstyle \binom{n}{k}}
\;\;\;\;\mbox{for all $k=1,2,\ldots ,n-2$.}
\end{equation}
The $n-2$ upper bounds in (\ref{sandwich}) are attained simultaneously
if and only if $T$ has precisely one leaf; the $n-2$ lower bounds in
(\ref{sandwich}) are attained simultaneously if and only if $T$ has
precisely $n-1$ leaves, in which case the root has $n-1$ children.
How high can these lower bounds be raised if an upper bound smaller
than $n-1$ is imposed on the number of children of every node of $T$?
This is the question we are going to answer. Its creative
interpretation has appeared in \cite{F03}.\\

Harary and Schwenk \cite{HS} call a tree (an unrooted one) where
removal of all vertices of degree one produces a path a {\em
caterpillar\/}. Abusing this usage a little, we will call a
caterpillar any rooted tree where removal of all leaves produces a
rooted tree with precisely one leaf. By a {\em full caterpillar of
degree $d$\/}, we will mean a caterpillar where each internal node,
except possibly the lowest one, has precisely $d$ children. 
\begin{theorem}\label{main}
Let $n$ and $d$ be positive integers such that $d<n$; let $T$ be any
rooted tree on $n$ nodes where each internal node has at most $d$
children and let $T^\ast$ be the full caterpillar of degree $d$ on $n$
nodes. Then $c(T,k)\ge c(T^\ast,k)$ for all $k=1,2,\ldots ,n$.
\end{theorem}
In the special case where $d=2$, inequalities $c(T,k)\ge c(T^\ast,k)$
with (essentially) $k=1,2,\ldots ,1+\lfloor \log_2 n\rfloor$ were
proved in \cite{F06} by an argument different from the argument given
below.\\

Given rooted trees $T,T'$ on $n$ nodes, we will write $T\succ T'$ to
mean that $c(T,k)\ge c(T',k)$ for all $k=1,2,\ldots ,n$ and that that
$c(T,k)>c(T',k)$ for at least one of these values of $k$. With this
notation, a refinement of Theorem~\ref{main} can be stated as follows.
\begin{theorem}\label{main2}
Let $n$ and $d$ be positive integers such that $d<n$ and let $T^\ast$
be the full caterpillar of degree $d$ on $n$ nodes.  If $T$ is a
rooted tree on $n$ nodes such that each internal node of $T$ has at
most $d$ children, then $T\succ T^\ast$ or else $T=T^\ast$.
\end{theorem}

Our proof of Theorem~\ref{main2} relies on two ways of altering a
rooted tree $T$ so that the resulting tree succeeds $T$ in the partial
order $\succ\,$. We shall describe these alterations in terms of the
{\em parent function\/} of a rooted tree that assigns to each node $z$
of the tree its parent $p(z)$ --- except when $z$ is the root, in
which case $p(z)$ is undefined.

\begin{lemma}\label{lift}
Let $T$ be a rooted tree defined by parent function $p$. Let $x$ and
$y$ be nodes of $T$ such that $x$ is not the root and $y$ is a
proper ancestor of $p(x)$. Let $T'$ be the rooted tree defined by
parent function $p'$ such that
\[
p'(z)=
\left\{
\begin{array}{ll}
p(z) & \mbox{if $z\ne x$,}\\
y & \mbox{if $z=x$.}
\end{array}
\right.
\]
Then $T\succ T'$.
\end{lemma}

\noindent{\bf Proof.} If $z$ is a leaf of $T$, then $z$ is a leaf of $T'$
and every node on the path from the root to $z$ in $T'$ lies on the
path from the root to $z$ in $T$. It follows that every transversal in
$T'$ is a transversal in $T$, and so $c(T,k)\ge c(T',k)$ for all $k$.
To see that $c(T,k)>c(T',k)$ for at least one $k$, consider the set
that consists of $p(x)$ and all leaves of $T$ that are not descendants
of $p(x)$: this set is a transversal in $T$ but not in $T'$.
\mbox{\hspace{0.5cm}}\hfill {\sc qed}

\begin{lemma}\label{shed}
Let $T$ be a rooted tree defined by parent function $p$. Let $x$ and
$y$ be nodes of $T$ such that $x$ is not the root, $x$ is not a leaf,
and $y$ is a leaf which is a proper descendant of a sibling of
$x$. Let $T'$ be the rooted tree defined by parent function $p'$ such
that
\[
p'(z)=
\left\{
\begin{array}{ll}
p(z) & \mbox{if $p(z)\ne x$,}\\
y & \mbox{if $p(z)=x$.}
\end{array}
\right.
\]
Then $T\succ T'$.
\end{lemma}

\noindent{\bf Proof.}
Given any set $S'$ of nodes in $T'$, define 
\[
f(S')=
\left\{
\begin{array}{ll}
S' & \mbox{if $S'$ meets the path from the root to $y$,}\\
S'-\{x\}\cup\{y\} & \mbox{otherwise.}
\end{array}
\right.
\]
If $S'$ is a transversal in $T'$, then $f(S')$ is a transversal in $T$
and $|f(S')|=|S'|$. If $f(S_1)=f(S_2)$ and $S_1\ne S_2$, then at least
one of $S_1,S_2$ avoids the path from the root to $x$, and so it is
not a transversal in $T'$. It follows that $c(T,k)\ge c(T',k)$ for all
$k$.  To see that $c(T,k)>c(T',k)$ for at least one $k$, consider the
set that consists of $p(y)$ and all leaves of $T$ that are not
descendants of $p(y)$: this set $S$ is a transversal in $T$, but there
is no transversal $S'$ in $T'$ such that $f(S')=S$.
\mbox{\hspace{0.5cm}}\hfill {\sc qed}\\

\noindent{\bf Proof of Theorem~\ref{main2}.} Consider any rooted tree $T$
on $n$ nodes such that each internal node of $T$ has at most $d$
children. Assuming that there is no rooted tree $T'$ on $n$ nodes such
that each internal node of $T'$ has at most $d$ children and such
that $T\succ T'$, we shall prove that $T$ is the full caterpillar of
degree $d$. Lemma~\ref{shed} guarantees that no two
internal nodes of $T$ are siblings, which means that $T$ is a
caterpillar; in turn, Lemma~\ref{lift} guarantees that each internal
node of $T$, except possibly the lowest one, has precisely $d$
children. \mbox{\hspace{0.5cm}}\hfill {\sc qed}\\

Imposing an upper bound on the number of children of every node is a
way of staying clear of the tree that attains simultaneously the $n-2$
lower bounds in (\ref{sandwich}), one where all children of the root
are leaves. Another way to stay clear of this tree is to impose an
upper bound on the number of leaves. There is a corresponding analogue
of Theorem~\ref{main2} and this analogue also follows directly from
our two lemmas.

\begin{theorem}
Let $n$ and $m$ be positive integers such that $m<n$ and let $T^\ast$
be the caterpillar on $n$ nodes, where the root has $m$ children and
every node other than the root has at most one child.  If $T$ is a
rooted tree with $n$ nodes and at most $m$ leaves, then $T\succ
T^\ast$ or else $T=T^\ast$.
\end{theorem}

\noindent{\bf Proof.} Consider any rooted tree $T$ with $n$ nodes and
at most $m$ leaves.  Assuming that there is no rooted tree $T'$ with
$n$ nodes and at most $m$ leaves such that $T\succ T'$, we shall prove
that $T=T^\ast$. Lemma~\ref{shed} guarantees that no two internal
nodes of $T$ are siblings, which means that $T$ is a caterpillar; in
turn, Lemma~\ref{lift} guarantees that no internal node of $T$ other
than the root has two or more children. \mbox{\hspace{0.5cm}}\hfill
{\sc qed}\\

\begin{center}
{\bf Acknowledgments}\\
\end{center}
We thank David Avis for inviting Jonathan David Farley to visit McGill
and we thank Jonathan David Farley for telling us about his nice
conjecture (= Theorem~\ref{main} with $d=2$) during this visit.


\begin{thebibliography}{99}

\bibitem{F03} J.D.~Farley, Breaking Al Qaeda cells: a
mathematical analysis of counterterrorism operations (a guide for risk
assessment and decision making), {\em Studies in Conflict and
Terrorism\/} {\bf 26} (2003), 399--411.

\bibitem{F06} J.D.~Farley, Toward a mathematical theory of counterterrorism. How to
build the perfect terrorist cell, I. Manuscript, 2006.


\bibitem{HS} F.~Harary and A.J.~Schwenk, Trees with hamiltonian
  square, {\em Mathematika\/} {\bf 18} (1971), 138--140.


\end{thebibliography}
\end{document}